\documentclass[preprint,9pt]{elsarticle}
\usepackage{amssymb}
\usepackage{amsmath}
\usepackage{amsthm}
\usepackage{enumerate}

\newtheorem{thrm}{Theorem}[section]

\newtheorem{exam}[thrm]{Example}
\newtheorem{cor}[thrm]{Corollary}

\theoremstyle{definition}

\newtheorem{remark}[thrm]{Remark}

\journal{...}
\input{tcilatex}

\begin{document}

\begin{frontmatter}


\cortext[cor1]{Corresponding author (+903562521616-3087)}

\title{Determinant and Permanent of Hessenberg Matrix and Generalized Lucas Polynomials}


\author[rvt]{Kenan Kaygisiz\corref{cor1}}\ead{kenan.kaygisiz@gop.edu.tr}
\author[rvt]{Adem Sahin} \ead{adem.sahin@gop.edu.tr}

\address[rvt]{Department of Mathematics, Faculty of Arts and Sciences,
Gaziosmanpa\c{s}a University, 60250 Tokat, Turkey}

\begin{abstract}
In this paper, we give some determinantal and permanental
representations of Generalized Lucas Polynomials by using various
Hessenberg matrices, which are general form of determinantal
and permanental representations of ordinary Lucas and Perrin
sequences. Then we show, under what conditions that the
determinants of the Hessenberg matrix becomes its permanents.
\end{abstract}
\begin{keyword}
Lucas numbers, Perrin numbers, generalized Lucas polynomials,
generalized Perrin polynomials, Hessenberg matrix, determinant and
permanent.\MSC[2010]{Primary 11B37, 15A15, Secondary 15A51}
\end{keyword}

\end{frontmatter}



\section{Introduction}

\bigskip Fibonacci numbers, Lucas numbers and Perrin numbers are%
\begin{equation*}
F_{n}=F_{n-1}+F_{n-2}\text{ for\ }n>2\text{ \ and }F_{1}=F_{2}=1,
\end{equation*}%
\begin{equation*}
L_{n}=L_{n-1}+L_{n-2}\text{ \ for }n>1\text{ \ and }L_{0}=2,\text{ }L_{1}=1,
\end{equation*}%
\begin{equation*}
R_{n}=R_{n-2}+R_{n-3}\text{ \ for }n>3\text{ \ and }R_{0}=3,\text{ }R_{1}=0,%
\text{ }R_{2}=2
\end{equation*}%
respectively.

There are large amount of studies on these sequences. In addition,
generalization of these sequences have been studied by many researchers.

Miles [10] defined generalized order-$k$ Fibonacci numbers(GO$k$F) as,%
\begin{equation}
f_{k,n}=\sum\limits_{j=1}^{k}f_{k,n-j}\
\end{equation}%
for $n>k\geq 2$, with boundary conditions: $f_{k,1}=f_{k,2}=f_{k,3}=\cdots
=f_{k,k-2}=0$ and $f_{k,k-1}=f_{k,k}=1.$\newline

\bigskip

Er [2] defined $k$ sequences of generalized order-$k$ Fibonacci numbers ($k$%
SO$k$F) as; for $n>0,$ $1\leq i\leq k$%
\begin{equation}
f_{k,n}^{\text{ }i}=\sum\limits_{j=1}^{k}c_{j}f_{k,n-j}^{\text{ }i}\ \
\label{mil}
\end{equation}%
with boundary conditions for $1-k\leq n\leq 0,$

\begin{equation*}
f_{k,n}^{\text{ }i}=\left\{
\begin{array}{l}
1\text{ \ \ \ \ \ if \ }i=1-n, \\
0\text{ \ \ \ \ \ otherwise,}%
\end{array}%
\right.
\end{equation*}%
where $c_{j}$ $(1\leq j\leq k)$ are constant coefficients, $f_{k,n}^{\text{ }%
i}$ is the $n$-th term of $i$-th sequence of order $k$ generalization. For $%
c_{j}=1$, $k$-th sequence of this generalization involves the Miles
generalization(\ref{mil}) for $i=k,$ i.e.%
\begin{equation}
f_{k,n}^{k}=f_{k,k+n-2}.  \label{fib}
\end{equation}

Kili\c{c} and Ta\c{c}c\i \lbrack 3] defined $k$ sequences of generalized
order-$k$ Pell numbers ($k$SO$k$P) as; for $n>0,$ $1\leq i\leq k$%
\begin{equation}
p_{k,n}^{\text{ }i}=2p_{k,n-1}^{\text{ }i}+p_{k,n-2}^{\text{ }i}+\cdots +\
p_{k,n-k}^{\text{ }i}\   \label{pell}
\end{equation}%
with initial conditions for $1-k\leq n\leq 0,$

\begin{equation*}
p_{k,n}^{\text{ }i}=\left\{
\begin{array}{l}
1\text{ \ \ \ \ \ if \ }i=1-n, \\
0\text{ \ \ \ \ \ otherwise,}%
\end{array}%
\right.
\end{equation*}%
where $p_{k,n}^{\text{ }i}$ is the $n$-th term of $i$-th sequence of order $%
k $ generalization.

\bigskip

MacHenry [7] defined generalized Fibonacci polynomials $(F_{k,n}(t))$, Lucas
polynomials $(G_{k,n}(t))$, where $t_{i}$ $(1\leq i\leq k)$ are constant
coefficients of the core polynomial%
\begin{equation}
P(x;t_{1},t_{2},\ldots ,t_{k})=x^{k}-t_{1}x^{k-1}-\cdots -t_{k},
\label{core}
\end{equation}%
which is denoted by the vector
\begin{equation}
t=(t_{1},t_{2},\ldots ,t_{k}).  \label{vkt}
\end{equation}%
$F_{k,n}(t)$ is defined inductively by
\begin{eqnarray}
F_{k,n}(t) &=&0,\text{ }n<1  \label{fippol} \\
F_{k,1}(t) &=&1  \notag \\
F_{k,2}(t) &=&t_{1}  \notag \\
F_{k,n+1}(t) &=&t_{1}F_{k,n}(t)+\cdots +t_{k}F_{k,n-k+1}(t).  \notag
\end{eqnarray}

$G_{k,n}(t_{1},t_{2},\ldots ,t_{k})$ is defined by%
\begin{eqnarray}
G_{k,n}(t) &=&0,\text{ }n<0  \label{lucpol} \\
G_{k,0}(t) &=&\text{ }k  \notag \\
G_{k,1}(t) &=&t_{1}\text{ }  \notag \\
G_{k,n+1}(t) &=&t_{1}G_{k,n}(t)+\cdots +t_{k}G_{k,n-k+1}(t).  \notag
\end{eqnarray}

Moreover in [8,9], MacHenry obtained some properties of these polynomials.
In [9] MacHenry gave a relation between generalized Fibonacci and Lucas
polynomials as; \
\begin{equation*}
G_{k,0}(t)=k,G_{k,n}(t)=F_{k,n+1}(t)+\sum_{j=1}^{k-1}jt_{j+1}F_{k,n-j}(t).
\end{equation*}

Equivalently this relation can be written as;%
\begin{equation*}
G_{k,0}(t)=k,\text{ }G_{k,n}(t)=\sum_{j=1}^{k}jt_{j}F_{k,n-j+1}(t).
\end{equation*}

\bigskip Kaygisiz and \c{S}ahin [4] definied generalized Perrin polynomials
by using generalized Lucas Polynomials. For $k\geq 3;$
\begin{eqnarray}
R_{k,0}(t) &=&k  \notag \\
R_{k,1}(t) &=&0  \notag \\
R_{k,2}(t) &=&2t_{2}  \notag \\
R_{k,3}(t) &=&t_{2}R_{k,1}(t)+3t_{3}  \label{perpol} \\
R_{k,4}(t) &=&t_{2}R_{k,2}(t)+t_{3}R_{k,1}(t)+4t_{4}  \notag \\
&&\vdots  \notag \\
R_{k,k-1}(t) &=&t_{2}R_{k,k-3}(t)+\cdots +t_{k-1}R_{k,1}(t)+kt_{k}  \notag
\end{eqnarray}%
and for $n\geq k$,%
\begin{equation*}
R_{k,n}(t)=\sum\limits_{i=2}^{k}t_{i}R_{k,n-i}(t).
\end{equation*}

\begin{remark}
Let $f_{k,k+n-2},$ $f_{k,n}^{\text{ }i},$ $p_{k,n}^{\text{ }i}$, $%
F_{k,n}(t), $ $G_{k,n}(t)$ and $R_{k,n}(t)$ be GO$k$F(\ref{mil}), $k$SO$k$F(%
\ref{fib}), $k$SO$k$P(\ref{pell}), generalized Fibonacci polynomials(\ref%
{fippol}), generalized Lucas polynomials(\ref{lucpol}) and generalized
Perrin polynomials(\ref{perpol}) respectively. Then \newline
$i)$ substituting $c_{i}=t_{i}$ for $1\leq i\leq k$ in (\ref{fib}) and (\ref%
{fippol}), we obtain the equality
\begin{equation*}
f_{k,n-1}^{\text{ }1}=F_{k,n}(t),
\end{equation*}%
$ii)$ substituting $t_{1}=2$ and $t_{i}=1$ for $2\leq i\leq k$ in (\ref%
{fippol}), we obtain the equality%
\begin{equation*}
p_{k,n-1}^{\text{ }1}=F_{k,n}(t),
\end{equation*}%
$iii)$ substituting $t_{1}=0$ in (\ref{lucpol}), we obtain the equality
\begin{equation*}
R_{k,n}(t)=G_{k,n}(t),
\end{equation*}%
$iv)$ substituting $t_{i}=1$ in (\ref{fippol}), we obtain the equality
\begin{equation*}
f_{k,k+n-2}=F_{k,n}(t),
\end{equation*}%
$v)$ substituting $t_{i}=1$ and $k=2$ in (\ref{lucpol}), we obtain the
equality%
\begin{equation*}
L_{n}=G_{k,n}(t),
\end{equation*}%
$vi)$ substituting $t_{1}=0$ and $t_{i}=1$ for $2\leq i\leq k$ and $k=3$ in (%
\ref{lucpol}), we obtain the equality%
\begin{equation*}
R_{n}=G_{k,n}(t).
\end{equation*}
\end{remark}

\bigskip

This Remark shows that $F_{k,n}(t)$ and $G_{k,n}(t)$ are general form of all
sequences mentioned above. Therefore, any result obtained from the
polynomial $F_{k,n}(t)$ and $G_{k,n}(t)$ are valid for other sequences.

Many researchers studied on determinantal and permanental representations of
$k$ sequences of generalized order-$k$ Fibonacci and Lucas numbers. For
example, Minc [11] defined an $n\times n$ (0,1)-matrix $F(n,k),$ and showed
that the permanents of $F(n,k)$ is equal to the generalized order-$k$
Fibonacci numbers.

In [5] and [6] the authors defined two (0,1)-matrices and showed that the
permanents of these matrices are the generalized Fibonacci and Lucas
numbers. \"{O}cal [12] gave some determinantal and permanental
representations of $k$-generalized Fibonacci and Lucas numbers and obtained
Binet's formulas for these sequences. Y\i lmaz and Bozkurt [15] derived some
relationships between Pell and Perrin sequences, and permanents and
determinants of a type of Hessenberg matrices. In [13] and [14] the authors
give some relation between determinant and permanent.

In this paper we give some determinantal and permanental representations of
Generalized Lucas Polynomials by using various Hessenberg matrices. In
addition we show under what conditions that the determinants of the
Hessenberg matrix becomes its permanent.

\section{The determinantal representations}

\bigskip

An $n\times n$ matrix $A_{n}=(a_{ij})$ is called lower Hessenberg matrix if $%
a_{ij}=0$ when $j-i>1$ i.e.,%
\begin{equation}
A_{n}=\left[
\begin{array}{ccccc}
a_{11} & a_{12} & 0 & \cdots & 0 \\
a_{21} & a_{22} & a_{23} & \cdots & 0 \\
a_{31} & a_{32} & a_{33} & \cdots & 0 \\
\vdots & \vdots & \vdots &  & \vdots \\
a_{n-1,1} & a_{n-1,2} & a_{n-1,3} & \cdots & a_{n-1,n} \\
a_{n,1} & a_{n,2} & a_{n,3} & \cdots & a_{n,n}%
\end{array}%
\right]  \label{an}
\end{equation}

\begin{thrm}
\label{cahil}$\bigskip \lbrack 1]$ $A_{n}$ be the $n\times n$ lower
Hessenberg matrix for all $n\geq 1$ and define $\det (A_{0})=1,$ then,%
\begin{equation*}
\det (A_{1})=a_{11}
\end{equation*}%
and for $n\geq 2$%
\begin{equation}
\det (A_{n})=a_{n,n}\det
(A_{n-1})+\sum\limits_{r=1}^{n-1}((-1)^{n-r}a_{n,r}\prod%
\limits_{j=r}^{n-1}a_{j,j+1}\det (A_{r-1})).  \label{det}
\end{equation}
\end{thrm}

\bigskip

\begin{thrm}
\label{t1}Let $k\geq 2$ be an integer, $G_{k,n}(t)$ be the generalized Lucas
Polynomials and $C_{k,n}=(c_{rs})$ be an $n\times n$ Hessenberg matrix, where%
\begin{equation*}
c_{rs}=\left\{
\begin{array}{l}
i^{\left\vert r-s\right\vert }.\frac{t_{r-s+1}}{t_{2}^{(r-s)}}\text{ \ \ \ \
\ \ \ \ \ \ \ \ \ \ \ \ \ \ \ \ \ \ \ \ if \ }s\neq 1\text{\ and }-1\leq
r-s<k\text{ }, \\
i^{\left\vert r-s\right\vert }.\frac{t_{r-s+1}}{t_{2}^{(r-s)}}.(r-s+1)\text{
\ \ \ \ \ \ \ if \ }s=1\text{\ and }-1\leq r-s<k\text{ },\text{\ \ \ } \\
0\text{ \ \ \ \ \ \ \ \ \ \ \ \ \ \ \ \ \ \ \ \ \ \ \ \ \ \ \ \ \ \ \ \ \ \
\ \ \ \ otherwise\ \ \ \ \ \ \ \ \ \ \ \ \ \ \ \ \ \ \ \ \ \ \ \ \ \ \ \ \ \
\ \ \ \ }%
\end{array}%
\right.
\end{equation*}%
i.e.,%
\begin{equation*}
Q_{k,n}=\left[
\begin{array}{cccccc}
t_{1} & it_{2} & 0 & 0 & \cdots  & 0 \\
2i & t_{1} & it_{2} & 0 & \cdots  & 0 \\
3i^{2}\frac{t_{3}}{t_{2}^{2}} & i & t_{1} & it_{2} & \cdots  & 0 \\
\vdots  & \vdots  & \vdots  & \vdots  &  & \vdots  \\
ki^{k-1}\frac{t_{k}}{t_{2}^{k-1}} & i^{k-2}\frac{t_{k-1}}{t_{2}^{k-2}} &
i^{k-3}\frac{t_{k-2}}{t_{2}^{k-3}} & i^{k-4}\frac{t_{k-3}}{t_{2}^{k-4}} &
\cdots  & 0 \\
0 & i^{k-1}\frac{t_{k}}{t_{2}^{k-1}} & i^{k-2}\frac{t_{k-1}}{t_{2}^{k-2}} &
i^{k-3}\frac{t_{k-2}}{t_{2}^{k-3}} & \cdots  & 0 \\
\vdots  & \vdots  & \vdots  & \vdots  & \ddots  & it_{2} \\
0 & 0 & 0 & \cdots  & i & t_{1}%
\end{array}%
\right] .
\end{equation*}%
Then%
\begin{equation*}
\det (C_{k,n})=G_{k,n}(t)
\end{equation*}%
where $t_{0}=1$ and $i=\sqrt{-1}.$
\end{thrm}

\bigskip

\begin{proof}
\bigskip Proof is by mathematical induction on $n$. The result is true for $%
n=1$ by hypothesis.

Assume that it is true for all positive integers less than or equal
to $m,$ that is $\det (C_{k,m})=G_{k,m+1}(t).$ Using Theorem
\ref{cahil} we have
\begin{eqnarray*}
\det (C_{k,m+1}) &=&q_{m+1,m+1}\det
(C_{k,m})+\sum\limits_{r=1}^{m}\left(
(-1)^{m+1-r}q_{m+1,r}\prod\limits_{j=r}^{m}q_{j,j+1}\det
(C_{k,r-1})\right)
\\
&=&t_{1}\det (C_{k,m})+\sum\limits_{r=1}^{m-k+1}\left(
(-1)^{m+1-r}q_{m+1,r}\prod\limits_{j=r}^{m}q_{j,j+1}\det
(C_{k,r-1})\right)
\\
&&+\sum\limits_{r=m-k+2}^{m}\left(
(-1)^{m+1-r}q_{m+1,r}\prod\limits_{j=r}^{m}q_{j,j+1}\det
(C_{k,r-1})\right)
\\
&=&t_{1}\det (C_{k,m})+\sum\limits_{r=m-k+2}^{m}\left(
(-1)^{m+1-r}q_{m+1,r}\prod\limits_{j=r}^{m}q_{j,j+1}\det
(C_{k,r-1})\right)
\\
&=&t_{1}\det (C_{k,m})+\sum\limits_{r=m-k+2}^{m}\left( (-1)^{m+1-r}.i^{m+1-r}%
\frac{t_{m-r+2}}{t_{2}^{(m-r+1)}}\prod\limits_{j=r}^{m}it_{2}\det
(C_{k,r-1})\right)  \\
&=&t_{1}\det (C_{k,m}) \\
&&+\sum\limits_{r=m-k+2}^{m}\left( (-1)^{m+1-r}.i^{m+1-r}\frac{t_{m-r+2}}{%
t_{2}^{(m-r+1)}}.i^{m+1-r}.t_{2}^{(m-r+1)}\det (Q_{k,r-1})\right)  \\
&=&t_{1}\det (Q_{k,m})+\sum\limits_{r=m-k+2}^{m}\left(
(-1)^{m+1-r}.i^{m+1-r}t_{m-r+2}.i^{m+1-r}.\det (Q_{k,r-1})\right)  \\
&=&t_{1}\det (Q_{k,m})+\sum\limits_{r=m-k+2}^{m}t_{m-r+2}\det (Q_{k,r-1}) \\
&=&t_{1}\det (Q_{k,m})+t_{2}\det (Q_{k,m-1})+\cdots +t_{k}\det
(Q_{k,m-(k-1)}).
\end{eqnarray*}%
From the hypothesis and the definition of generalized Lucas
polynomials we
obtain%
\begin{equation*}
\det (C_{k,m+1})=t_{1}G_{k,m}(t)+t_{2}G_{k,m-1}(t)+\cdots
+t_{k}G_{k,m-(k-1)}(t)=G_{k,m+1}(t).
\end{equation*}%
Therefore, the result is true for all possitive integers.
\end{proof}

\bigskip

\begin{exam}
We obtain $6$-th Generalized Lucas polynomials for $k=5,$ i.e. $G_{5,6}(t),$
by using Theorem \ref{t1}.%
\begin{eqnarray*}
\det (C_{5,6}) &=&\det \left[
\begin{array}{cccccc}
t_{1} & -it_{2} & 0 & 0 & 0 & 0 \\
2i & t_{1} & -it_{2} & 0 & 0 & 0 \\
3\frac{-t_{3}}{t_{2}^{2}} & i & t_{1} & -it_{2} & 0 & 0 \\
4\frac{-it_{4}}{t_{2}^{3}} & \frac{-t_{3}}{t_{2}^{2}} & i & t_{1} & -it_{2}
& 0 \\
5\frac{t_{5}}{t_{2}^{4}} & \frac{-it_{4}}{t_{2}^{3}} & \frac{-t_{3}}{%
t_{2}^{2}} & i & t_{1} & -it_{2} \\
0 & \frac{t_{5}}{t_{2}^{4}} & \frac{-it_{4}}{t_{2}^{3}} & \frac{-t_{3}}{%
t_{2}^{2}} & i & t_{1}%
\end{array}%
\right]  \\
&=&6t_{1}t_{5}+6t_{2}t_{4}+12t_{1}t_{2}t_{3}+2t_{2}^{3}+3t_{3}^{2}+%
\allowbreak
t_{1}^{6}+6t_{1}^{2}t_{4}+6t_{1}^{3}t_{3}+6t_{1}^{4}t_{2}+9t_{1}^{2}t_{2}^{2}\allowbreak
\\
&=&G_{5,6}(t).
\end{eqnarray*}
\end{exam}

\bigskip

\begin{thrm}
\label{t2}\bigskip Let $k\geq 2$ be an integer$,$ $G_{k,n}$ be the
generalized Lucas Polynomial and $B_{k,n}=(b_{ij})$ be an $n\times n$ lower
Hessenberg matrix such that%
\begin{equation*}
b_{ij}=\left\{
\begin{array}{l}
-t_{2}\text{ \ \ \ \ \ \ \ \ \ \ \ \ \ \ \ \ \ \ \ \ \ \ \ \ \ \ if \ \ \ }%
j=i+1, \\
\frac{t_{i-j+1}}{t_{2}^{(i-j)}}\text{\ \ \ \ \ \ \ \ \ \ \ \ \ \ \ \ \ \ \ \
\ \ \ \ \ if\ \ }i\neq 1\text{\ \ and }0\leq i-j<k\text{,} \\
\frac{t_{i-j+1}}{t_{2}^{(i-j)}}.(i-j+1)\text{\ \ \ \ \ \ \ \ if\ \ \ }i=1%
\text{\ \ and }0\leq i-j<k\text{,} \\
0\text{\ \ \ \ \ \ \ \ \ \ \ \ \ \ \ \ \ \ \ \ \ \ \ \ \ \ \ \ \ \ otherwise}%
\end{array}%
\right.
\end{equation*}%
i.e.,%
\begin{equation*}
B_{k,n}=\left[
\begin{array}{cccccc}
t_{1} & -t_{2} & 0 & 0 & \cdots  & 0 \\
2 & t_{1} & -t_{2} & 0 & \cdots  & 0 \\
3\frac{t_{3}}{t_{2}^{2}} & 1 & t_{1} & -t_{2} & \cdots  & 0 \\
\vdots  & \vdots  & \vdots  & \vdots  &  & \vdots  \\
k\frac{t_{k}}{t_{2}^{k-1}} & \frac{t_{k-1}}{t_{2}^{k-2}} & \frac{t_{k-2}}{%
t_{2}^{k-3}} & \ddots  & \cdots  & 0 \\
0 & \frac{t_{k}}{t_{2}^{k-1}} & \frac{t_{k-1}}{t_{2}^{k-2}} & \frac{t_{k-2}}{%
t_{2}^{k-3}} & \cdots  & 0 \\
\vdots  & \vdots  & \vdots  & \vdots  & \ddots  & -t_{2} \\
0 & 0 & 0 & \cdots  & \cdots  & t_{1}%
\end{array}%
\right] .
\end{equation*}%
Then%
\begin{equation*}
\det (B_{k,n})=G_{k,n}(t).
\end{equation*}%
where $t_{0}=1.$
\end{thrm}

\bigskip

\begin{proof}
\bigskip \bigskip Proof is similar to the proof of Theorem \ref{t1} using
Theorem \ref{cahil}.
\end{proof}

\begin{exam}
\bigskip We obtain $5$-th generalized Lucas polynomial for $k=4$, i.e. $%
G_{4,5}(t),$ by using Theorem \ref{t2}.%
\begin{eqnarray*}
B_{4,5} &=&\left[
\begin{array}{ccccc}
t_{1} & -t_{2} & 0 & 0 & 0 \\
2 & t_{1} & -t_{2} & 0 & 0 \\
3\frac{t_{3}}{t_{2}^{2}} & 1 & t_{1} & -t_{2} & 0 \\
4\frac{t_{4}}{t_{2}^{3}} & \frac{t_{3}}{t_{2}^{2}} & 1 & t_{1} & -t_{2} \\
0 & \frac{t_{4}}{t_{2}^{3}} & \frac{t_{3}}{t_{2}^{2}} & 1 & t_{1}%
\end{array}%
\right]  \\
&=&5t_{1}t_{4}+5t_{2}t_{3}+t_{1}^{5}+5t_{1}t_{2}^{2}+5t_{1}^{2}t_{3}+%
\allowbreak 5t_{1}^{3}t_{2} \\
&=&G_{4,5}(t).
\end{eqnarray*}
\end{exam}

\begin{cor}
\bigskip \bigskip\ If we rewrite Theorem \ref{t1} and Theorem \ref{t2} for $%
t_{i}=1$ and $k=2,$ we obtain%
\begin{equation*}
\det (C_{k,n})=L_{n}^{\text{ }}
\end{equation*}%
and%
\begin{equation*}
\det (B_{k,n})=L_{n}^{\text{ }}
\end{equation*}%
respectively, where $L_{n}^{\text{ }}$ are the ordinary Lucas numbers.
\end{cor}

\bigskip

\begin{cor}
\bigskip \bigskip If we rewrite Theorem \ref{t1} and Theorem \ref{t2} for $%
t_{1}=0$ for $1\leq i\leq k,$ we obtain%
\begin{equation*}
\det (C_{k,n})=R_{k,n}(t)
\end{equation*}%
and%
\begin{equation*}
\det (B_{k,n})=R_{k,n}(t)
\end{equation*}%
respectively, where $R_{k,n}(t)$ are the generalized Perrin polynomials.
\end{cor}

\bigskip

\begin{cor}
If we rewrite Theorem \ref{t1} and Theorem \ref{t2} for $t_{1}=0$ and $%
t_{i}=1$ for $2\leq i\leq k$ and $k=3$ we obtain%
\begin{equation*}
\det (C_{k,n})=R_{n}
\end{equation*}%
and%
\begin{equation*}
\det (B_{k,n})=R_{n}
\end{equation*}%
respectively, where $R_{n}$ are the ordinary Perrin numbers.
\end{cor}

\bigskip

\section{The permanent representations}

\bigskip

Let $A=(a_{i,j})$ be an $n\times n$ square matrix over a ring R. The
permanent of $A$ is defined by%
\begin{equation*}
\text{per}(A)=\sum\limits_{\sigma \in
S_{n}}\prod\limits_{i=1}^{n}a_{i,\sigma (i)}
\end{equation*}%
where $S_{n}$ denotes the symmetric group on $n$ letters.

\bigskip

\begin{thrm}
\label{ocal}$\left[ 10\right] $Let $A_{n}$ be an $n\times n$ lower
Hessenberg matrix for all $n\geq 1$ and define per$(A_{0})=1.$ Then,%
\begin{equation*}
\text{per}(A_{1})=a_{11}
\end{equation*}%
and for $n\geq 2$%
\begin{equation}
\text{per}(A_{n})=a_{n,n}\text{per}(A_{n-1})+\sum\limits_{r=1}^{n-1}(a_{n,r}%
\prod\limits_{j=r}^{n-1}a_{j,j+1}\text{per}(A_{r-1})).  \label{per}
\end{equation}
\end{thrm}

\bigskip

\begin{thrm}
\label{t3}\bigskip \bigskip Let $k\geq 2$ be an integer, $G_{k,n}(t)$ be the
generalized Lucas Polynomials and $H_{k,n}=(h_{rs})$ be an $n\times n$ lower
Hessenberg matrix such that%
\begin{equation*}
h_{rs}=\left\{
\begin{array}{l}
i^{(r-s)}.\frac{t_{r-s+1}}{t_{2}^{(r-s)}}\text{ \ \ \ \ \ \ \ \ \ \ \ \ \ \
\ \ \ \ \ \ \ \ \ \ \ if \ }s\neq 1\text{\ and}-1\leq r-s<k, \\
i^{(r-s)}.\frac{t_{r-s+1}}{t_{2}^{(r-s)}}.(r-s+1)\text{ \ \ \ \ \ \ \ if \ }%
s=1\text{\ and }-1\leq r-s<k\text{ },\text{\ \ \ \ \ \ \ \ \ \ \ \ \ \ \ \ \
\ \ \ \ \ \ \ \ \ \ \ \ \ \ } \\
0\text{ \ \ \ \ \ \ \ \ \ \ \ \ \ \ \ \ \ \ \ \ \ \ \ \ \ \ \ \ \ \ \ \ \ \
\ \ \ \ \ \ otherwise\ \ \ \ \ \ \ \ \ \ \ \ \ \ \ \ \ \ \ \ \ \ \ \ \ \ \ \
\ \ \ \ \ \ }%
\end{array}%
\right.
\end{equation*}%
i.e.,%
\begin{equation}
H_{k,n}=\left[
\begin{array}{cccccc}
t_{1} & -it_{2} & 0 & 0 & \cdots  & 0 \\
2i & t_{1} & -it_{2} & 0 & \cdots  & 0 \\
3i^{2}\frac{t_{3}}{t_{2}^{2}} & i & t_{1} & -it_{2} & \cdots  & 0 \\
\vdots  & \vdots  & \vdots  & \vdots  &  & \vdots  \\
ki^{k-1}\frac{t_{k}}{t_{2}^{k-1}} & i^{k-2}\frac{t_{k-1}}{t_{2}^{k-2}} &
i^{k-3}\frac{t_{k-2}}{t_{2}^{k-3}} & i^{k-4}\frac{t_{k-3}}{t_{2}^{k-4}} &
\cdots  & 0 \\
0 & i^{k-1}\frac{t_{k}}{t_{2}^{k-1}} & i^{k-2}\frac{t_{k-1}}{t_{2}^{k-2}} &
i^{k-3}\frac{t_{k-2}}{t_{2}^{k-3}} & \cdots  & 0 \\
\vdots  & \vdots  & \vdots  & \vdots  & \ddots  &  \\
0 & 0 & 0 & \cdots  & \cdots  & t_{1}%
\end{array}%
\right] .
\end{equation}%
Then%
\begin{equation*}
\text{per}(H_{k,n})=G_{k,n}(t)
\end{equation*}%
where $t_{0}=1$ and $i=\sqrt{-1}.$
\end{thrm}

\bigskip

\begin{proof}
\bigskip Proof is similar to the proof of Theorem \ref{t1} by using Theorem
\ref{ocal}.
\end{proof}

\bigskip

\begin{exam}
\bigskip We obtain $3$-th Generalized Lucas polynomials for $k=4$, by using
Theorem \ref{t3}%
\begin{eqnarray*}
\det (H_{4,3}) &=&\det \left[
\begin{array}{ccc}
t_{1} & -it_{2} & 0 \\
2i & t_{1} & -it_{2} \\
3\frac{-t_{3}}{t_{2}} & i & t_{1}%
\end{array}%
\right]  \\
&=&3t_{3}+3t_{1}t_{2}+t_{1}^{3}.
\end{eqnarray*}
\end{exam}

\bigskip

\begin{thrm}
\label{t4}Let $k\geq 2$ be an integer$,$ $G_{k,n}(t)$ be the generalized
Lucas Polynomials and $L_{k,n}=(l_{ij})$ be an $n\times n$ lower Hessenberg
matrix such that%
\begin{equation*}
l_{ij}=\left\{
\begin{array}{l}
\begin{array}{l}
\frac{t_{i-j+1}}{t_{2}^{(i-j)}}\text{\ \ \ \ \ \ \ \ \ \ \ \ \ \ \ \ \ \ \ \
\ \ \ \ \ if\ \ }j\neq 1\text{\ \ and }0\leq i-j<k\text{,} \\
\frac{t_{i-j+1}}{t_{2}^{(i-j)}}.(i-j+1)\text{\ \ \ \ \ \ \ \ if\ \ \ }j=1%
\text{\ \ and }0\leq i-j<k\text{,}%
\end{array}
\\
0\text{ \ \ \ \ \ \ \ \ \ \ \ \ \ otherwise\ \ \ \ \ \ \ \ \ \ \ \ \ \ \ \ \
\ \ \ \ \ \ \ \ \ \ \ \ \ \ \ \ \ }%
\end{array}%
\right.
\end{equation*}%
i.e.,%
\begin{equation*}
L_{k,n}=\left[
\begin{array}{cccccc}
t_{1} & t_{2} & 0 & 0 & \cdots  & 0 \\
2 & t_{1} & t_{2} & 0 & \cdots  & 0 \\
3\frac{t_{3}}{t_{2}^{2}} & 1 & t_{1} & t_{2} & \cdots  & 0 \\
\vdots  & \vdots  & \vdots  & \vdots  &  & \vdots  \\
k\frac{t_{k}}{t_{2}^{k-1}} & \frac{t_{k-1}}{t_{2}^{k-2}} & \frac{t_{k-2}}{%
t_{2}^{k-3}} & \frac{t_{k-3}}{t_{2}^{k-4}} & \cdots  & 0 \\
0 & \frac{t_{k}}{t_{2}^{k-1}} & \frac{t_{k-1}}{t_{2}^{k-2}} & \frac{t_{k-2}}{%
t_{2}^{k-3}} & \cdots  & 0 \\
& \vdots  & \vdots  & \vdots  & \ddots  &  \\
0 & 0 & 0 & \cdots  & \cdots  & t_{1}%
\end{array}%
\right]
\end{equation*}%
where $t_{0}=1.$ Then%
\begin{equation*}
\text{per}(L_{k,n})=G_{k,n}(t).
\end{equation*}
\end{thrm}

\bigskip

\begin{proof}
\bigskip Proof of the theorem is similar to the proof of Theorem \ref{t1}
using Theorem \ref{ocal}.
\end{proof}

\bigskip

\begin{cor}
If we rewrite Theorem \ref{t3} and Theorem \ref{t4} for $t_{i}=1$ $(1\leq
i\leq k)$ and $k=2$ $,$ we obtain%
\begin{equation*}
\text{per}(H_{k,n})=L_{n}
\end{equation*}%
and%
\begin{equation*}
\text{per}(L_{k,n})=L_{n}
\end{equation*}%
respectively, where $L_{n}$ are the ordinary Lucas numbers.
\end{cor}

\bigskip

\begin{cor}
\bigskip If we rewrite Theorem \ref{t3} and Theorem \ref{t4} for $t_{1}=0$
and $t_{i}=1$ $(2\leq i\leq k),$ we obtain
\begin{equation*}
\text{per}(H_{k,n})=R_{k,n}(t)
\end{equation*}%
and%
\begin{equation*}
\text{per}(L_{k,n})=R_{k,n}(t)
\end{equation*}%
respectively, where $R_{k,n}(t)$ are the generalized Perrin polynomials.
\end{cor}

\bigskip

\begin{cor}
\bigskip If we rewrite Theorem \ref{t3} and Theorem \ref{t4} for $t_{1}=0$, $%
t_{i}=1$ $(2\leq i\leq k)$ and $k=3$ we obtain
\begin{equation*}
\text{per}(H_{k,n})=R_{n}
\end{equation*}%
and%
\begin{equation*}
\text{per}(L_{k,n})=R_{n}
\end{equation*}%
respectively, where $R_{n}$ are the ordinary Perrin numbers.
\end{cor}

\subsection{Determinat and Permanent of a Hessenberg Matrix}

In this section we give a relation between the determinant and the permanent
of a Hessenberg matrix.

\begin{thrm}
Let $A_{n}$\bigskip\ be the Hessenberg matrix in (\ref{an}) and $%
B_{n}=(b_{ij})$ be an $n\times n$ Hessenberg matrix such that%
\begin{equation*}
b_{ij}=\left\{
\begin{array}{l}
0\text{\ \ \ \ \ \ \ \ \ \ \ \ if \ }j-i>1, \\
-a_{ij}\text{ \ \ \ \ \ \ \ if \ \ }j-i=1\text{ },\text{\ \ \ \ \ \ \ \ \ \
\ \ \ \ \ \ \ \ \ \ \ \ \ \ \ \ \ \ \ \ \ } \\
a_{ij}\text{ \ \ \ \ \ \ \ \ \ otherwise.\ \ \ \ \ \ \ \ \ \ \ \ \ \ \ \ \ \
\ \ \ \ \ \ \ \ \ \ \ \ \ \ }%
\end{array}%
\right.
\end{equation*}%
i.e.,%
\begin{equation*}
B_{n}=\left[
\begin{array}{ccccc}
a_{11} & -a_{12} & 0 & \cdots & 0 \\
a_{21} & a_{22} & -a_{23} & \cdots & 0 \\
a_{31} & a_{32} & a_{33} & \cdots & 0 \\
\vdots & \vdots & \vdots &  & \vdots \\
a_{n-1,1} & a_{n-1,2} & a_{n-1,3} & \cdots & -a_{n-1,n} \\
a_{n,1} & a_{n,2} & a_{n,3} & \cdots & a_{n,n}%
\end{array}%
\right].
\end{equation*}%
Then%
\begin{equation*}
\det B_{n}=\text{per}A_{n}
\end{equation*}%
or%
\begin{equation*}
\det A_{n}=\text{per}B_{n}.
\end{equation*}
\end{thrm}

\bigskip

\begin{proof}
\bigskip We know from (\ref{det})
\begin{equation*}
\det (A_{1})=a_{11}
\end{equation*}%
and for $n\geq 2$%
\begin{equation*}
\det (A_{n})=a_{n,n}\det
(A_{n-1})+\sum\limits_{r=1}^{n-1}((-1)^{n-r}a_{n,r}\prod%
\limits_{j=r}^{n-1}a_{j,j+1}\det (A_{r-1}))
\end{equation*}%
and from (\ref{per})
\begin{equation*}
\text{per}(A_{1})=a_{11}
\end{equation*}%
and for $n\geq 2$%
\begin{equation*}
\text{per}(A_{n})=a_{n,n}\text{per}(A_{n-1})+\sum\limits_{r=1}^{n-1}(a_{n,r}%
\prod\limits_{j=r}^{n-1}a_{j,j+1}\text{per}(A_{r-1})).
\end{equation*}%
Using mathematical induction on $n$, we prove this theorem by using (\ref{det}) and
(\ref{per}). The result is true
for $n=1$ by hypothesis.

Assume that it is true for all positive integers less than or equal
to $m,$ namely $\det B_{m}=$per$A_{m}.$ For $n\geq 2$
\begin{eqnarray*}
\det (B_{m+1}) &=&a_{m+1,m+1}\det
(B_{m})+\sum\limits_{r=1}^{m}((-1)^{m+1-r}a_{m+1,r}\prod%
\limits_{j=r}^{m}b_{j,j+1}\det (B_{r-1})) \\
&=&a_{m+1,m+1}\text{per}(A_{m})+\sum\limits_{r=1}^{m}((-1)^{m+1-r}a_{m+1,r}%
\prod\limits_{j=r}^{m}(-a_{j,j+1})\text{per}(A_{r-1})) \\
&=&a_{m+1,m+1}\text{per}(A_{m})+\sum%
\limits_{r=1}^{m}((-1)^{m+1-r}a_{m+1,r}(-1)^{m+1-r}\prod%
\limits_{j=r}^{m}a_{j,j+1}\text{per}(A_{r-1})) \\
&=&a_{m+1,m+1}\text{per}(A_{m})+\sum\limits_{r=1}^{m}(a_{m+1,r}\prod%
\limits_{j=r}^{m+1}a_{j,j+1}\text{per}(A_{r-1})) \\
&=&\text{per}(A_{m+1}).
\end{eqnarray*}%
Therefore, the result is true for all possitive integers.
\end{proof}

\bigskip 

\bigskip

\end{document}